\documentclass[12pt, twoside]{article}
\usepackage{amsmath,amsthm,amssymb}
\usepackage{times}
\usepackage{enumerate}

\def\titlerunning#1{\gdef\titrun{#1}}
\makeatletter
\def\author#1{\gdef\autrun{\def\and{\unskip, }#1}\gdef\@author{#1}}
\def\address#1{{\def\and{\\\hspace*{18pt}}\renewcommand{\thefootnote}{}%
\footnote {#1}}%
\markboth{\autrun}{\titrun}}
\makeatother
\def\email#1{e-mail: #1}
\def\subjclass#1{{\renewcommand{\thefootnote}{}%
\footnote{\emph{Mathematics Subject Classification (2010):} #1}}}
\def\keywords#1{\par\medskip
\noindent\textbf{Keywords.} #1}

\newtheorem{thm}{Theorem}[section]
\newtheorem{cor}[thm]{Corollary}
\newtheorem{lem}[thm]{Lemma}

\theoremstyle{definition}
\newtheorem{defin}[thm]{Definition}
\newtheorem{rem}[thm]{Remark}

\numberwithin{equation}{section}

\begin{document}


\baselineskip=17pt


\titlerunning{Toda System}

\title{Local profile of fully bubbling solutions to SU(n+1) Toda Systems}

\author{Chang-Shou Lin
\and
Juncheng Wei
\and
Lei Zhang}

\date{}

\maketitle

\address{Lin:Taida Institute of Mathematical Sciences
and Center for Advanced Study in Theoretical Sciences,
        National Taiwan University,
         Taipei 106, Taiwan ; \email{cslin@math.ntu.edu.tw}
\and
Wei: Department of Mathematics,
        Chinese University of Hong Kong,
        Shatin, Hong Kong  and
        Department of Mathematics, University of British Columbia,
        Vancouver, V6P 1Z2; \email{wei@math.cuhk.edu.hk}
\and
Zhang: Department of Mathematics,
        University of Florida,
        358 Little Hall P.O.Box 118105,
        Gainesville FL 32611-8105, \email{leizhang@ufl.edu}}



\begin{abstract}
In this article we prove that for locally defined singular $SU(n+1)$ Toda systems in $\mathbb R^2$, the profile of fully bubbling solutions near the singular source can be accurately approximated by global solutions. The main ingredients of our new approach are the classification theorem of Lin-Wei-Ye \cite{lin-wei-ye} and the non-degeneracy of the linearized Toda system \cite{lin-wei-ye}, which make us overcome the difficulties that come from lack of symmetry and the singular source.
\end{abstract}

\section{Introduction}

Let $(M,g)$ be a compact Riemann surface and $\Delta$ the Beltrami-Laplacian operator of the metric $g$, and $K$ the Gauss curvature. The $SU(n+1)$ Toda system is the following nonlinear PDE

\begin{equation}\label{toda1}
\Delta u_i+\sum_{j=1}^n a_{ij} h_j e^{u_j} -K(x)=4\pi \sum_j \gamma_{ij} \delta_{q_j}, \quad 1\le i\le n,
\end{equation}
where $h_i$ ($i=1,...,n$) are positive smooth functions on $M$, $\delta_q$ stands for the Dirac measure at $q\in M$, and
$A=(a_{ij})$ is the Cartan matrix given by
$$
 A=\left(\begin{array}{ccccc}
2 & -1 & 0 & ... & 0\\
-1 & 2 & -1 & ... & 0\\
0 & -1 & 2 &  & 0\\
\vdots & \vdots & & \vdots\\
0 & \ldots & -1 & 2 & -1\\
0 & \ldots &  & -1 & 2
\end{array}
\right).
$$

Toda system (\ref{toda1}) has aroused a lot of attention in recent years because of its close connection to many different fields of mathematics and physics.
 For $n=1$, (\ref{toda1}) is reduced to the Gauss curvature equation in two dimensional surfaces. Without the singular source and $M=\mathbb S^2$, it is the well known Nirenberg problem. In general it is related to the existence of metric of positive constant curvature with conic singularities (\cite{chenli1, chenli2, troy1, troy2}). For the past three decades, equation (\ref{toda1}) with $n=1$ has been extensively studied (see \cite{chang},  \cite{mfdegree}, \cite{licmp}  for example). For the general $n$ and $h_i\equiv 1$ ($i=1,..,n$) equation (\ref{toda1}) is connected with holomorphic curves of $M$ into $\mathbb CP^n$ via the classical infinitesimal Pl\"ucker formulae, see \cite{h-griffith}. This geometric connection is very important because from it, it has been found out that equation (\ref{toda1}) with $h_i\equiv 1$ is an integrable system ( see \cite{chern},\cite{guest}, for example).
Recently by using this connection, Lin-Wei-Ye \cite{lin-wei-ye} are able to completely classify all the entire solutions of (\ref{toda1}) in $\mathbb R^2$ with one singular source and finite energy.

In mathematical physics, equation (\ref{toda1}) has also played an important role in Chern-Simons gauge theory. For example, in the relativistic $SU(n+1)$ Chern-Simons model proposed by physicists ( see \cite{jackiw}) for $n=1$ and \cite{dunne2} for $n>1$), in order to explain the physics of high temperature super-conductivity, (\ref{toda1}) governs the limiting equations as physical parameters tend to $0$. For the past twenty years, the connections of (\ref{toda1}) with $n=1$ and the Chern-Simons-Higgs equation have been explored extensively. See \cite{nolasco1} and \cite{ lin-yan1}. However, for $n\ge 2$ only very few works are devoted to this direction of research. See \cite{ao-lin-wei}, \cite{ lin-yan2} and \cite{nolasco2}.
For recent development of equation (\ref{toda1}) and related subjects, we refer the readers to
\cite{bjm, bm,  jostlinwang,jostwang,lin-wei-ye,lwz1,lin-wei-zhao,linzhang1,linzhang2,linzhang3, mal-ruiz1, mal-ruiz2, yang2} and the reference therein.

One of the fundamental issues concerning (\ref{toda1}) is to study the bubbling phenomenon, which could lead to establishing a priori bound of solutions of (\ref{toda1}). For the case $n=1$, the bubbling phenomenon has been studied thoroughly for the past twenty years. Basically there are two kinds of bubbling behaviors of solutions near its blowup points. One is called ``simple blowup", which means the bubbling profile could be well controlled locally by entire bubbling solutions in $\mathbb R^2$. For the case without singular sources, this was proved by Y. Y. Li \cite{licmp}, applying the method of moving planes. If there is a singular source $4\pi\gamma \delta_0$ on the right hand side of the equation, this was proved by Bartolucci-Chen-Lin-Tarantello \cite{bartolucci4} for $\gamma\not\in \mathbb N$, and recently by Kuo-Lin \cite{kuo-lin} if $\gamma\in \mathbb N$, who use potential analysis and Pohozaev identity. On the other hand, the non-simple blowup could occur at $\gamma\in \mathbb N$ only. The sharp profile of the non-simple blowup has recently been proved in \cite{kuo-lin}. The study of the bubbling phenomenon is important not only for deriving a priori bounds, but also for providing a lot of important geometric information near blowup points, see \cite{chen-lin-1,chen-lin-2,lin-yan2}.

For $n\ge 2$, (\ref{toda1}) is an elliptic system. It is expected that the behavior of bubbling solutions is more complicated than the case $n=1$. One major difficulty comes from the partial blown-up phenomenon, that is, after a suitable scaling, the solutions converge to a smaller system. To understand the partial blown-up phenomenon, we have to first study the fully blown-up behavior, and to obtain accurate description of this class of bubbling solutions.  When $n=2$ and (\ref{toda1}) has no singular sources, the bubbling behavior of fully bubbling solutions has been studied by Jost-Lin-Wang \cite{jostlinwang} and
Lin-Wei-Zhao \cite{lin-wei-zhao}.  In \cite{jostlinwang} it is proved that any sequence of fully bubbling solutions is a simple blowup at any blowup point. The proof in \cite{jostlinwang} uses deep application of holonomy theory, which is a very effective generalization of Pohozaev identity. Unfortunately their holonomy method cannot be extended to cover the case with singular sources.  The purpose of this article is to extend their results to  any $n\geq 2$ and to include (\ref{toda1}) with singular sources. Before stating our main results, we set up our problem first. Since this is a local problem, for simplicity we consider

\begin{equation}\label{su3s}
 \Delta u_i^k+\sum_{j=1}^n a_{ij} h_j^ke^{u_j^k}=4\pi \gamma_i\delta_0, \quad B_1\subset \mathbb R^2
 \end{equation}
 where $B_1$ is the unit ball. We shall use $B_r$ to denote the ball centered at origin with radius $r$.

 For $u^k=(u_1^k,..,u_n^k)$, $h^k=(h_1^k,..,h_n^k)$ and $\gamma_i$ ($i=1,..,n$) we
 assume the usual assumptions:
 \begin{eqnarray*}
 (H):\quad &(i):& \frac 1C\le h_i^k\le C, \quad \|h_i^k\|_{C^2(B_1)}\le C, \quad h_i^k(0)=1, \quad i=1,..,n\\
 &(ii):& \gamma_i>-1,\quad i=1,..,n \\
 &(iii):& \int_{B_1}h_i^ke^{u_i^k}\le C,\quad i=1,..,n, \mbox{ $C$ is independent of $k$}. \\
 &(iv):& |u_i^k(x)-u_i^k(y)|\le C, \quad \mbox{for all } x,y\in \partial B_1, \quad i=1,..,n. \\
 &(v):& \max_{K\subset\subset B_1\setminus \{0\}} u_i^k\le C, \mbox{ and $0$ is the only blowup point}.
 \end{eqnarray*}

If $(u_1^k,...,u_n^k)$ is a global solution of (\ref{toda1}) in $M$, it is easy to see that all assumptions of $(H)$ are satisfied.
We also note that the assumption $(iv)$ in $(H)$ is necessary for our analysis, without it Chen \cite{chenxx} proved that even for $n=1$
the blowup solutions can be very complicated near their blowup points. The assumption $h_i^k(0)=1$ in $(i)$ is just for convenience.

Let
\begin{equation}\label{epk}
-2\log \epsilon_k=\max_{x\in B_1,i=1,..,n}(\frac{\tilde u_i^k(x)}{1+\gamma_i}),\quad \mbox{ where }
\tilde u_i^k(x)=u_i^k(x)-2\gamma_i \log |x|,
\end{equation}
and
\begin{equation}\label{tvik}
\tilde v_i^k(y)=\tilde u_i^k(\epsilon_k y)+2(1+\gamma_i)\log \epsilon_k, \quad i=1,..,n
\end{equation}
Then clearly $\tilde v_i^k$ satisfies
\begin{equation}\label{vik}
\Delta \tilde v_i^k(y)+\sum_{j=1}^n a_{ij}|y|^{2\gamma_j}h_j^k(\epsilon_ky)e^{\tilde v_j^k}=0,\quad |y|\le \epsilon_k^{-1}.
\end{equation}
Our major assumption is $\tilde v^k=(\tilde v_1^k,...,\tilde v_n^k)$ converges to a $SU(n+1)$ Toda system uniformly over all compact subsets of $\mathbb R^2$:
\begin{defin}\label{fully-bu}
We say $u^k$ of (\ref{su3s}) is a fully bubbling sequence if $\tilde v^k$ converges in $C^{1,\alpha}_{loc}(\mathbb R^2)$ to $\tilde v=(\tilde v_1,..,\tilde v_n)$ that solves the following $SU(n+1)$ Toda system in $\mathbb R^2$:
\begin{eqnarray}\label{singv}
\Delta \tilde v_i+\sum_{j=1}^n a_{ij} |y|^{2\gamma_j}e^{\tilde v_j}=0,\quad \mathbb R^2,\quad i=1,..,n\\
\int_{\mathbb R^2}|y|^{2\gamma_i}e^{\tilde v_i}<\infty, \quad i=1,..,n. \nonumber
\end{eqnarray}
\end{defin}

The main purpose of this paper is to show that a fully bubbling sequence $u^k$ can be sharply approximated by a sequence of global solutions
$U^k=(U_1^k,..,U_n^k)$ of
\begin{equation}\label{globalUk}
\Delta U_i^k+\sum_{j=1}^n a_{ij} e^{U_j^k}=4\pi \gamma_i\delta_0, \quad \mbox{ in }\quad \mathbb R^2,\quad i=1,...,n.
\end{equation}
\begin{thm}\label{thm-main}
Let (H) hold, $u^k$ be a fully bubbling
sequence described in Definition \ref{fully-bu} and $\epsilon_k$ be defined in (\ref{epk}). Then there exists a sequence of global solutions $U^k=(U_1^k,..,U_n^k)$ of (\ref{globalUk})
such that for $|y|\le \epsilon_k^{-1}$ and $i=1,...,n$
\begin{align}\label{main-est}
&|u_i^k(\epsilon_k y)-U_i^k(\epsilon_k y)|\\
\le & \left \{\begin{array}{ll}
C(\sigma)\epsilon_k^{\sigma}(1+|y|)^{\sigma}, \,\, &\mbox{ if } \min\{\gamma_1,...,\gamma_n\}\le -\frac 34,\,\,
\sigma\in (0,\min\{2+2\gamma_1,.,,,2+2\gamma_n\})\\
C\epsilon_k (1+|y|),\,\, &\mbox{ if } \min\{ \gamma_1,.., \gamma_n \}>-\frac 34.
\end{array}
\right.  \nonumber
\end{align}
Moreover, there exists $C>0$
independent of $k$, such that
\begin{equation}\label{global-ap}
|\tilde U_i^k(\epsilon_k y)+2(1+\gamma_i)\log \epsilon_k+2(2+\gamma_i+\gamma_{n+1-i})\log (1+|y|)|\le C,
\end{equation}
for $|y|\le \epsilon_k^{-1}$ and $i=1,...,n$,
where
$\tilde U_i^k(x)=U_i^k(x)-2\gamma_i\log |x|$ is the regular part of $U_i^k$.
\end{thm}

The global solutions
\begin{equation}\label{tuik}
(\tilde U_1^k(\epsilon_ky)+2(1+\gamma_1)\log \epsilon_k,..,\tilde U_n^k(\epsilon_ky)+2(1+\gamma_n)\log \epsilon_k)
\end{equation}
in Theorem \ref{thm-main} are perturbations of $\tilde v=(\tilde v_1,...,\tilde v_n)$ in (\ref{singv}). In fact, the sequence in (\ref{tuik}) converges uniformly to $\tilde v$ over any fixed compact subset of $\mathbb R^2$. Thus Theorem \ref{thm-main} clearly leads to the following

\begin{cor} \label{cor1} Let $u^k$, $\epsilon_k$ be the same as in Theorem \ref{thm-main},  $\tilde v^k$ be defined by (\ref{tvik}). Then for $i=1,..,n$,
\begin{equation}\label{ap-cor}
|\tilde v_i^k(y)+2(2+\gamma_i+\gamma_{n+1-i})\log (1+|y|)|\le C, \quad \mbox{ for } |y|\le \epsilon_k^{-1}.
\end{equation}
\end{cor}

\begin{rem} The estimate in (\ref{ap-cor}) holds trivially over any fixed compact subset of $\mathbb R^2$. So the strength of Corollary \ref{cor1} lies on the fact that the estimate is over $|y|\le \epsilon_k^{-1}$. Such type of estimate was first established by Li \cite{licmp} for single Liouville equations.
\end{rem}

Estimates similar to (\ref{main-est}) and (\ref{ap-cor}) can be found in \cite{licmp,chen-lin-1,bartolucci4,zhangcmp,zhangccm} for single Liouville equations and \cite{jostlinwang,lin-wei-zhao} for Toda systems. The proof of Theorem \ref{thm-main} is almost entirely different from all the approaches in these works. For example the estimates for single Liouville equations use ODE theory, which is based on the symmetry of global solutions. Lin-Wei-Zhao \cite{lin-wei-zhao}'s sharp estimates are tailored for regular $SU(3)$ Toda systems because they need to differentiate
blowup solutions at blowup points twice (which cannot be expected when the singular source exists) and a lot of algebraic computation to fix Cauchy data of blowup solutions. For general singular
$SU(n+1)$ Toda system, first, ODE method cannot be used because global solutions may not have any symmetry. Second, fixing Cauchy data of blowup solutions at a blowup point is impossible, because in addition to the differentiation issue mentioned before, the amount of algebraic computation required to fix Cauchy data depends on $n^2+2n$ parameters and is extremely complicated if $n$ is large.
 Our approach is purely based on PDE methods and the essential part relies on and important classification theorem of Lin-Wei-Ye \cite{lin-wei-ye} for global $SU(n+1)$ Toda system and the non-degeneracy property of the corresponding linearized system. The key point is to choose a sequence of global solutions as approximating solutions. On one hand these global solutions all tend to the limit system (\ref{singv}), which means all the $n^2+2n$ families of parameters corresponding to these global solutions have limit. On the other hand, one component of the approximating global solutions is very close to the same component of blowup solutions at $n^2+2n$ carefully chosen points. The closeness in one component leads to the closeness in other components as well.

Theorem \ref{thm-main} is an extension of previous works. For example, if
$n=2$ and $\gamma_i=0 (i=1,2)$, Corollary \ref{cor1} was proved by Jost-Lin-Wang \cite{jostlinwang}. It is easy to see that Theorem \ref{thm-main} is stronger than Corollary \ref{cor1} even for this special case. Lin-Wei-Zhao proved (\ref{main-est}) for $n=2$ and $\gamma_i=0 (i=1,2)$ but Theorem \ref{thm-main} also holds when then number of equation is greater than $2$ and the singular source at $0$ exists.

For some applications such as constructing blowup solutions, more refined estimates than those in Theorem \ref{thm-main} are needed. For $SU(3)$ Toda systems with no singularity, Lin-Wei-Zhao \cite{lin-wei-zhao} obtained more delicate estimates for this case based on Corollary \ref{cor1}.

The organization of the article is as follows. In section two we list some facts on the $SU(n+1)$ Toda system and the non-degeneracy of the linearized system. The proof of Theorem \ref{thm-main} is in section three. One key point in the proof of Theorem \ref{thm-main} is to determine $n^2+2n$ points in $\mathbb R^2$ in a specific way. Since this part is somewhat elaborate and elementary, we put it separately in section four.

\section{Some facts on the linearized $SU(n+1)$ system}

First we list some facts on the entire solutions of $SU(n+1)$ Toda systems with singularities. For more details see \cite{lin-wei-ye}.
Let $u=(u_1,..,u_n)$ solve
\begin{equation}\label{efor-u}
\left\{\begin{array}{ll}
\Delta u_i+\sum_{j=1}^n a_{ij} e^{u_j}=4\pi \gamma_i \delta_0, \quad \mathbb R^2,\quad i=1,..,n \\
\int_{\mathbb R^2} e^{u_i}<\infty
\end{array}
\right.
\end{equation}
where $A=(a_{ij})_{n\times n}$ is the Cartan matrix and $\gamma_i>-1$. Then let
$$u^i=\sum_{j=1}^n a^{ij}u_j,\quad i=1,...,n $$
where $(a^{ij})_{n\times n}=A^{-1}$.
Clearly $(u^1,...,u^n)$ satisfies
$$\Delta u^i+e^{\sum_{j=1}^n a_{ij} u^j}=4\pi \gamma^i \delta_0, \quad \mbox{where} \,\, \gamma^i=\sum_{j=1}^n a^{ij} \gamma_j, \quad i=1,...,n.\quad  $$

The classification theorem of Lin-Wei-Ye (\cite{lin-wei-ye}) asserts
\begin{equation}\label{13may28e1}
e^{-u^1}=|z|^{-2\gamma^1}(\lambda_0+\sum_{i=1}^n \lambda_i |P_i(z)|^2)
\end{equation}
where for
$$\mu_i=1+\gamma_i,\quad i=1,...,n $$
\begin{equation}\label{piz}
P_i(z)=z^{\mu_1+...+\mu_i}+\sum_{j=0}^{i-1}c_{ij} z^{\mu_1+...+\mu_j}, \quad i=1,...,n
\end{equation}
 $c_{ij}$ ($j<i$) are complex numbers and $\lambda_i>0$ ($0\le i\le n$) satisfies
\begin{equation}\label{nlam}
\lambda_0...\lambda_n=2^{-n(n+1)}\Pi_{1\le i\le j\le n}(\sum_{k=i}^j\mu_k)^{-2}.
\end{equation}
Furthermore if $\mu_{j+1}+...+\mu_i\not \in \mathbb N$ for some $j<i$, $c_{ij}=0$.
Let
$$\tilde u^1=u^1-2\gamma^1 \log |z|, $$
then
\begin{equation}\label{tu1}
\tilde u^1=-\log (\lambda_0+\sum_{i=1}^n \lambda_i |P_i(z)|^2).
\end{equation}

The following lemma classifies the solutions of the linearized system under a mild growth condition at infinity:
\begin{lem}\label{luniq}
Let $\Phi_1,...,\Phi_n$ solve the linearized $SU(n+1)$ Toda system:
\begin{equation}\label{linear-sys}
\Delta \Phi_i+e^{u_i}(\sum_{j=1}^n a_{ij} \Phi_j)=0, \quad \mbox{in}\quad \mathbb R^2,\quad i=1,...,n
\end{equation}
where $u$ solves (\ref{efor-u}).
If
\begin{equation}\label{13may28e5}
|\Phi_i(x)|\le C(1+|x|)^{\sigma}, \quad x\in \mathbb R^2
\end{equation}
for $\sigma\in (0, \min \{1, 2\mu_1,...,2\mu_n\})$, then
\begin{equation}\label{13may28e2}
e^{-u^1}\Phi_1(z)=\sum_{k=0}^n m_{kk}|z|^{2\beta_k}
+2\sum_{k=1}^{n-1}|z|^{\beta_k}\sum_{l=k+1}^n |z|^{\beta_l}Re(m_{kl}e^{-i(\mu_{k+1}+..+\mu_l)\theta})
\end{equation}
where $\theta=arg(z)$,
\begin{equation}
\label{betadefinition}
\beta_0=-\gamma^1,\quad \beta_i=\gamma^i-\gamma^{i+1}+i,\quad \beta_n=\gamma^n+n,
\end{equation}
$m_{kk}\in \mathbb R$ for $k=0,..,n$, $m_{kl}\in \mathbb C$ for $k<l$. Obviously $m_{kl}=0$ if $\mu_{k+1}+..+\mu_l\not \in \mathbb N$.
\end{lem}

\noindent{\bf Proof of Lemma \ref{luniq}:}  This lemma is proved in \cite{lin-wei-ye} when all $\Phi_i$ are bounded functions. Here we mention the minor modifications when a mild growth condition in (\ref{13may28e5}) is assumed. Let
$$w_i(y)=-\frac 1{2\pi}\int_{\mathbb R^2} \log |y-\eta |e^{u_i(\eta)}(\sum_{j=1}^n a_{ij} \Phi_j(\eta))d\eta. $$
By (\ref{13may28e2}) and $e^{u_i(z)}=O(|z|^{-4-2\nu_{n+1-i}})$ we see that
$e^{u_i(z)}(\sum_{j=1}^n a_{ij} \Phi_j(z))=O(|z|^{-2-\delta})$ for some $\delta>0$ when $|z|$ is large. Thus $w_i(y)=O(\log |y|)$ for $|y|$ large.
From $\Delta (\Phi_i-w_i)=0$ in $\mathbb R^2$ and $|\Phi_i(z)-w_i(z)|\le O(|z|^{1-\delta})$ for some $\delta>0$ we have
$$\Phi_i=w_i+C. $$
Then using the integral representation of $\Phi_i$ we can further obtain $\nabla^k\Phi_i=O(|z|^{-k})$ as $|z|\to \infty$. Then the remaining part of the proof is the same as Lemma 6.1 of \cite{lin-wei-ye}. $\Box$

\medskip

From (\ref{betadefinition}) it is easy to verify that
\begin{equation}\label{betai}
\beta_i-\beta_{i-1}=\mu_i,\quad 1\le i\le n.
\end{equation}
Then we see that $\beta_i$ is increasing because $\mu_i=1+\gamma_i>0$.
Using (\ref{13may28e1}) and (\ref{betai}) in  (\ref{13may28e2}), we have
\begin{eqnarray}\label{13may28e3}
&&\Phi_1=\frac 1{\lambda_0+\sum_i \lambda_i |P_i(z)|^2}\bigg \{ \sum_{k=0}^n m_{kk}|z|^{2\beta_k+2\gamma^1}\\
&&+2\sum_{k=0}^{n-1}|z|^{\beta_k+\gamma^1}\sum_{l=k+1}^n |z|^{\beta_l+\gamma^1}Re(m_{kl}e^{-i(\mu_{k+1}+..+\mu_l)\theta})\bigg \}\nonumber\\
&&=\frac 1{\lambda_0+\sum_i \lambda_i |P_i(z)|^2}\bigg \{ \sum_{k=0}^n m_{kk} |z|^{2\mu_1+...+2\mu_k}+2\sum_{k=0}^{n-1}
|z|^{\mu_1+...+\mu_k}\nonumber\\
&&\bigg (\sum_{l=k+1}^n |z|^{\mu_1+...+\mu_l}Re(m_{kl}e^{-i(\mu_{k+1}+...+\mu_l)\theta})\bigg )\bigg \}. \nonumber
\end{eqnarray}
\medskip

\begin{lem}\label{lem1} $\displaystyle{\frac{m_{00}}{\lambda_0}+...+\frac{m_{nn}}{\lambda_n}=0}$.
\end{lem}

\noindent{\bf Proof of Lemma \ref{lem1}:} It is proved in \cite{lin-wei-ye} that the linearized system is non-degenerate, which means all solutions to (\ref{linear-sys}) are obtained by differentiating $n^2+2n$ parameters of $(u^1,...,u^n)$. In particular

\begin{equation}\label{phi-22}
\Phi_1=c_1\frac{\partial u^1}{\partial \lambda_1}+...+c_n\frac{\partial u^1}{\partial \lambda_n}+ c_{n+1}\frac{\partial u^1}{\partial c_{01}^{\mathbb R}}
+...c_{n^2+2n}\frac{\partial u^1}{\partial c_{n,n-1}^{\mathbb I}},
\end{equation}
where $c_{ij}^{\mathbb R}$ is the real part of $c_{ij}$, $c_{ij}^{\mathbb I}$ is the imaginary part.
Direct computation from (\ref{13may28e1}) and (\ref{nlam}) shows
$$\frac{\partial u^1}{\partial \lambda_k}=-\frac{|P_k|^2+\frac{\partial \lambda_0}{\partial \lambda_k}}{\lambda_0+\sum_{i=1}^n \lambda_i |P_i|^2}=-\frac{|P_k|^2-\frac{\lambda_0}{\lambda_k}}{\lambda_0+\sum_{i=1}^n \lambda_i |P_i|^2} $$
for $k=1,...,n$. Comparing (\ref{13may28e3}) and (\ref{phi-22}) we have
\begin{align*}
m_{kk}=-c_k, \quad k=1,...,n \\
m_{00}=\frac{c_1\lambda_0}{\lambda_1}+...\frac{c_n\lambda_0}{\lambda_n}.
\end{align*}
Then it is easy to see that
$$\frac{m_{00}}{\lambda_0}+\frac{m_{11}}{\lambda_1}+...+\frac{m_{nn}}{\lambda_n}=0. $$
  Lemma \ref{lem1} is established. $\Box$

\medskip

From Lemma \ref{lem1} we see that there are $n^2+2n$ unknowns in $\Phi_1$.  We write $\Phi_1$ as
\begin{eqnarray}\label{Phi1}
&&\Phi_1=\frac{1}{\lambda_0+\sum_{i=1}^n \lambda_i |P_i(z)|^2}\bigg \{
\sum_{k=1}^{n}m_{kk}|z|^{2\mu_1+..+2\mu_k}-\sum_{k=1}^n\frac{\lambda_0}{\lambda_k}m_{kk}\\
&&+2\sum_{k=0}^{n-1}|z|^{2\mu_1+..+2\mu_k}\sum_{l=k+1}^n Re(\bar m_{kl}z^{\mu_{k+1}+..+\mu_l})\bigg \}. \nonumber
\end{eqnarray}

\section{The Proof of Theorem \ref{thm-main}}

Recall that $\tilde v^k=(\tilde v_1^k,..,\tilde v_n^k)$ satisfies (\ref{vik}) and $\tilde v^k$ converges in $C^{1,\alpha}_{loc}(\mathbb R^2)$ to $\tilde v=(\tilde v_1,....,\tilde v_n)$ of (\ref{singv}). By the classification theorem of Lin-Wei-Ye \cite{lin-wei-ye}, there exists
$\Lambda=(\lambda_i, c_{ij}) (i=0,..,n,j<i)$ such that $\tilde v^1(z)$ is defined in (\ref{tu1}) where $\lambda_i$ and $P_i$ satisfy
(\ref{nlam}) and (\ref{piz}), respectively. To emphasize the dependence of $\Lambda$, we denote $\tilde v_i$ and $\tilde v^i$ as $\tilde v_i(z,\Lambda)$ and
$\tilde v^i(z,\Lambda)$, respectively.

The following matrix plays an important role in the argument below: For $p_1,...,p_{n^2+2n}\in \mathbb R^2$, set
\begin{equation}\label{thetamatrix}
{\bf M}= (\Theta(p_1),...,\Theta(p_{n^2+2n})).
\end{equation}
where
$$\Theta(p)=(\frac{\partial \tilde v^1}{\partial \lambda_0}(p),..,\frac{\partial \tilde v^1}{\partial \lambda_{n-1}}(p),
\frac{\partial \tilde v^1}{\partial c_{10}^{\mathbb R}}(p),...,\frac{\partial \tilde v^1}{\partial c_{n,n-1}^{\mathbb I}}(p))'. $$
where$( )'$ stands for transpose. 
In section four we shall show that by choosing $p_1,...,p_{n^2+2n}$ appropriately with respect to $\Lambda$ the matrix ${\bf M}$ is invertible.

Let $\tilde v^{i,k}=\sum_j a^{ij} \tilde v_j^k$, then $\tilde v^{i,k}$ converges uniformly to $\tilde v^i(\cdot, \Lambda)$ over any fixed compact subset of $\mathbb R^2$.
Since the difference between $\tilde v^{i,k}$ and $\tilde v^i(\cdot, \Lambda)$ is only $o(1)$, we need to find a sequence of global solutions that approximates better.
Suppose the sequence of global solutions is represented by $\Lambda_k:=(\lambda_i^k, c_{ij}^k)$: the regular part of the first component is
$$\tilde v^1(z,\Lambda_k)=-\log (\lambda_0^k+\sum_{i=1}^n \lambda_i^k |P_i^k(z)|^2) $$
with
$$P_i^k(z)=z^{\mu_1^k+...+\mu_i^k}+\sum_{j=0}^{i-1}c_{ij}^k z^{\mu_1^k+...+\mu_j^k}. $$
Other components $\tilde v^i(z,\Lambda_k)$ are determined by the equation
$$\Delta \tilde v^i(y,\Lambda_k)+|y|^{2\gamma_i}e^{\sum_j a_{ij} \tilde v^j(y,\Lambda_k)}=0,\quad \mbox{ in }\quad \mathbb R^2, \quad i=1,..,n. $$
Finally we set
\begin{equation}\label{global-vuk}
v^i(z,\Lambda_k)=\tilde v^i(z,\Lambda_k)(z)+2\gamma^{i}\log |z|,\,\,\mbox{where}\,\,\gamma^{i}=\sum_j a^{ij}\gamma_j,\,\, i=1,..,n.
\end{equation}
Then we claim that if
\begin{equation}\label{agree-p}
\tilde v^1(p_l,\Lambda_k)=\tilde v^{1,k}(p_l),\quad l=1,..,n^2+2n,
\end{equation}
we have
\begin{equation}\label{coincide}
\lambda_i^k\to \lambda_i,\,\, c_{ij}^k\to c_{ij}.
\end{equation}
Indeed, since $\tilde v^{1,k}(p_l)=\tilde v^1(p_l,\Lambda)+o(1)$ for $l=1,...,n^2+2n$, (\ref{coincide}) clearly follows from the invertibility of ${\bf M}$.
In other words there exists $\Lambda_k\to \Lambda$ such that (\ref{agree-p}) holds.

Let $v_i(\cdot,\Lambda_k)=\sum_ja_{ij}v^j(\cdot, \Lambda_k)$. Here we point out that
$$v_i(\cdot,\Lambda_k)=\tilde U_i^k(\epsilon_k \cdot)+2(1+\gamma_i)\log \epsilon_k, \quad i=1,...,n, $$
which is the global sequence in (\ref{tuik}) and the statement of Theorem \ref{thm-main}.

In order to obtain estimates (\ref{main-est}) we write (\ref{Phi1}) as
\begin{eqnarray}\label{Phi1a}
&&\Phi_1(z)(\lambda_0+\sum_i \lambda_i |P_i(z)|^2)\\
&=&\sum_{k=1}^{n}m_{kk}(|z|^{2\mu_1+..+2\mu_k}-\frac{\lambda_0}{\lambda_k})+2\sum_{k=0}^{n-1}\sum_{l=k+1}^n|z|^{2\mu_1+..+2\mu_k+\mu_{k+1}
+..+\mu_l}\nonumber \\
&& \quad (\cos((\mu_{k+1}+..+\mu_l)\theta) m_{kl}^1+\sin((\mu_{k+1}+..+\mu_l)\theta) m_{kl}^2). \nonumber \\
&=&{\bf X}\hat \Theta(z). \nonumber
\end{eqnarray}
where
$${\bf X}=(m_{11},...,m_{nn},m_{01}^1,...,m_{n-1,n}^2),\quad m_{kl}=m_{kl}^1+\sqrt{-1}m_{kl}^2. $$
So $\hat \Theta(z)$ is a column vector (so is $\Theta(p)$). Our choice of $p_1,...,p_{n^2+2n}$ ( explained in section four) also makes
$${\bf M_1}=(\hat \Theta(p_1),...,\hat \Theta(p_{n^2+2n})) $$
invertible.

Let $\Phi_i^k=\tilde v^{i,k}-\tilde v^i(\cdot,\Lambda_k)$. By (\ref{vik}) and the definition of $\tilde v^{i,k}$ we have
$$\Delta \tilde v^{i,k}+|y|^{2\gamma_i}h_i^k(\epsilon_ky)e^{\sum_j a_{ij} \tilde v^{j,k}(y)}=0, \quad |y|\le \epsilon_k^{-1}. $$
Hence the equation for $(\Phi_1^k,...,\Phi_n^k)$ can be written as
\begin{equation}\label{Phik}
\Delta \Phi_i^k(y)+|y|^{2\gamma_i}e^{\xi_i^k(y)}(\sum_j a_{ij} \Phi_j^k(y))=O(\epsilon_k|y|)|y|^{2\gamma_i}e^{\sum_j a_{ij} \tilde v^{j,k}}
\end{equation}
where, by the mean value theorem,
$$ e^{\xi_i^k}=\frac{e^{\sum_j a_{ij}\tilde v^{j,k}}-e^{\sum_j a_{ij}\tilde v^j(\cdot,\Lambda_k)}}{\sum_j a_{ij}(\tilde v^{j,k}-\tilde v^j(\cdot,\Lambda_k))}
=\int_0^1 e^{\sum_j a_{ij}(t\tilde v^{j,k}+(1-t)\tilde v^j(\cdot,\Lambda_k))}dt. $$
By Theorem 4.1 and Theorem 4.2 of \cite{lwz1}, $e^{\xi_i^k}$ converges uniformly to $e^{\tilde v_i(\cdot, \Lambda)}$ over all compact subsets of $\mathbb R^2$, moreover,
\begin{equation}\label{xik}
|y|^{2\gamma_i}e^{\xi_i^k(y)}=O(1+|y|)^{-4-2\gamma_{n+1-i}+o(1)}, \quad |y|\le \epsilon_k^{-1}.
\end{equation}
Also by Theorem 4.1 and Theorem 4.2 of \cite{lwz1} we can estimate the right hand side of (\ref{Phik}). Thus (\ref{Phik}) can be written as
\begin{equation}\label{Phik2}
\Delta \Phi_i^k+|y|^{2\gamma_i}e^{\xi_i^k(y)}(\sum_{j=1}^n a_{ij} \Phi_j^k)=\frac{O(\epsilon_k)}{(1+|y|)^{3+2\gamma_{n+1-i}}}, \,\,
\mbox{ in }\,\, |y|\le \epsilon_k^{-1}.
\end{equation}
It is immediate to observe that the oscillation of $\Phi_i^k$ on $\partial B_{\epsilon_k^{-1}}$ is finite. Thus for convenience we use the following functions to eliminate the oscillation of $\Phi_i^k$ on $\partial B_{\epsilon_k^{-1}}$:
$$\left\{\begin{array}{ll}
\Delta \psi_i^k=0, &\quad \mbox{ in }\quad B_{\epsilon_k^{-1}},\\
\psi_i^k=\Phi_i^k-\frac 1{2\pi \epsilon_k^{-1}}\int_{\partial B_{\epsilon_k^{-1}}} \Phi_i^k, &\quad \mbox{ on }\quad \partial B_{\epsilon_k^{-1}}.
\end{array}
\right.
$$
Standard estimate gives
\begin{equation}\label{psik}
|\psi_i^k(y)|\le C\epsilon_k |y|,\quad |y|\le \epsilon_k^{-1}.
\end{equation}
Let $\tilde \Phi_i^k=\Phi_i^k-\psi_i^k$, then by (\ref{Phik2}) and (\ref{psik}) we have
\begin{equation}\label{tphi}
\Delta \tilde \Phi_i^k+|y|^{2\gamma_i}e^{\xi_i^k(y)}(\sum_{j=1}^n a_{ij} \tilde \Phi_j^k)=\frac{O(\epsilon_k)}{(1+|y|)^{3+2\gamma_{n+1-i}}}, \,\,
\mbox{ in }\,\, |y|\le \epsilon_k^{-1}
\end{equation}
and it follows from (\ref{agree-p}) and (\ref{psik}) that
\begin{equation}\label{tphi1}
\tilde \Phi_1^k(p_l)=O(\epsilon_k), \quad l=1,...,n^2+2n.
\end{equation}
From here we consider two cases.

{\bf Case one: } $\min\{\gamma_1,...,\gamma_n\}\le -\frac 34. $

In this case we
set $$H_k=\max_{i}\max_{|y|\le \epsilon_k^{-1}}\frac{|\tilde \Phi_i^k(y)|}{(1+|y|)^{\sigma}\epsilon_k^{\sigma}} $$
for any fixed $\sigma\in (0,\min\{1,2\mu_1,..,2\mu_n\})$. Our goal is to show that $H_k$ is bounded. We prove this by contradiction. Suppose $H_k\to \infty$ and let $y_k$ be where the maximum is attained. Let
$$\hat \Phi_i^k(y)=\frac{\tilde \Phi_i^k(y)}{H_k (1+|y_k|)^{\sigma}\epsilon_k^{\sigma}}. $$
This definition immediately implies
\begin{equation}\label{onees}
|\hat \Phi_i^k(y)|=\frac{|\tilde \Phi_i^k(y)|}{H_k\epsilon_k^{\sigma}(1+|y|)^{\sigma}}\frac{(1+|y|)^{\sigma}}{(1+|y_k|)^{\sigma}}\le \frac{(1+|y|)^{\sigma}}{(1+|y_k|)^{\sigma}}.
\end{equation}
Next we write the equation for $(\hat \Phi_1^k,..,\hat \Phi_n^k)$ as
$$
\Delta \hat \Phi_i^k+|y|^{2\gamma_i}e^{\xi_i^k}(\sum_j a_{ij}\hat \Phi_j^k)=\frac{O(\epsilon_k^{1-\sigma})(1+|y|)^{-3-2\gamma_{n+1-i}}}{H_k(1+|y_k|)^{\sigma}},
$$
and we observe that
$\hat \Phi_i^k$  has no oscillation on $\partial B_{\epsilon_k^{-1}}$.

We first consider the case that along a subsequence, $y_k\to y^*$. In this case, $(\hat \Phi_1^k,..,\hat \Phi_n^k)$ converges to
$(\Phi_1,...,\Phi_n)$ that satisfies
\begin{equation}\label{13may29e1}
\left\{\begin{array}{ll}
\Delta \Phi_i+e^{v_i}\sum_j a_{ij}\Phi_j=0,\quad \mbox{ in } \quad \mathbb R^2, \quad i=1,...,n \\
\\
|\Phi_i(y)|\le C(1+|y|)^{\sigma},\quad i=,1..,n,\quad \sigma\in (0, \min\{1,2\mu_1,...,2\mu_n\}),\\
\\
\Phi_1(p_l)=0,\quad l=1,...,n^2+2n.
\end{array}
\right.
\end{equation}
where $v_i(y)=\tilde v_i(y)+2\gamma_i\log |y|$.
Note that the last equation in (\ref{13may29e1}) holds because of (\ref{tphi1}). From the first two equations of (\ref{13may29e1}) and
Lemma \ref{luniq} we have (\ref{13may28e2}).  Then by (\ref{Phi1a}) we have
$${\bf M} \hat \Theta(p_l)=0,\quad l=1,...,n^2+2n. $$
Since ${\bf M}$
is invertible, we have
$$m_{11}=...=m_{n,n}=m_{1,0}^1=...=m_{n,n-1}^2=0. $$
Thus $\Phi_1\equiv 0$, which means $\Phi_i\equiv 0$ for all $i$. This is a contradiction to
$|\Phi_i(y^*)|=1$ for some $i$.

The only remaining case we need to consider is when $y_k\to \infty$. To get a contradiction we evaluate
\begin{eqnarray}\label{13may14e1}
&&\hat \Phi_i^k(y_k)-\hat \Phi_i^k(0)\\
&=&\int_{B_{\epsilon_k^{-1}}}(G_k(y_k,\eta)-G_k(0,\eta))\bigg ( |\eta |^{2\gamma_i}e^{\xi_i^k(\eta)}(\sum_j a_{ij}\tilde \Phi_j^k(\eta))
\nonumber \\
&&\qquad \qquad +\frac{O(\epsilon_k^{1-\sigma})(1+|\eta |)^{-3-2\gamma_{n+1-i}}}{H_k(1+|y_k|)^{\sigma}}\bigg )d\eta
\nonumber
\end{eqnarray}
where $G_k$ is the Green's function on $B_{\epsilon_k^{-1}}$ with Dirichlet boundary condition.
To evaluate the right hand side of the term above we use (\ref{onees}),(\ref{xik})
and the following estimate of the Green's function (see \cite{linzhang3} for the proof) :

For $y\in \Omega_k:=B_{1/\epsilon_k}$,
let \begin{eqnarray*}
\Sigma_1&=&\{\eta \in \Omega_k;\quad |\eta |<|y|/2 \quad \}\\
\Sigma_2&=&\{\eta \in \Omega_k;\quad |y-\eta |<|y|/2 \quad \}\\
\Sigma_3&=&\Omega_k\setminus (\Sigma_1\cup \Sigma_2).
\end{eqnarray*}
Then for $|y|>2$,
\begin{equation}\label{1020e5}
|G_k(y,\eta)-G_k(0,\eta)|\le \left\{\begin{array}{ll}
C(\log |y|+|\log |\eta ||),\quad \eta\in \Sigma_1,\\
C(\log |y|+|\log |y-\eta ||),\quad \eta\in \Sigma_2,\\
C|y|/|\eta |,\quad \eta \in \Sigma_3.
\end{array}
\right.
\end{equation}

\medskip

Using (\ref{1020e5}) to estimate the right hand side of (\ref{13may14e1}) is standard. Here we just point out that we use (\ref{onees}) to estimate $\tilde \Phi_j^k(\eta)$ in the first term and it is essential to use $\epsilon_k^{1-\sigma}$ for the second term, as $\min\{2\mu_1,...,2\mu_n\}+\sigma$ may be less than or equal to $1$ in this case. At the end of these standard estimates we see that the right hand side of (\ref{13may14e1}) is $o(1)$. However we know $|\hat\Phi_i^k(y_k)|=1$ for some $i$ and it is easy to prove $|\hat\Phi_i^k(0)|\to 0$ by
exactly the same argument used in the proof of $y_k\to \infty$. Thus we obtain a contradiction and proved
 $$|\tilde \Phi_i^k(y)|\le C\epsilon_k^{\sigma}(1+|y|)^{\sigma}. $$

{ \bf Case two: } $\min\{\gamma_1,...,\gamma_n\}>-\frac 34. $

In this case we set
$$H_k=\max_{i}\max_{|y|\le \epsilon_k^{-1}}\frac{|\tilde \Phi_i^k(y)|}{(1+|y|)^{\sigma}\epsilon_k} $$
and
$$\hat \Phi_i^k(y)=\frac{\tilde \Phi_i^k(y)}{H_k(1+|y_k|)^{\sigma}}. $$
Here we choose $\sigma$ not only in $(0,\min\{1,2\mu_1,...,2\mu_n\})$, but also satisfy
\begin{equation}\label{13dec9e2}
\min\{2\mu_1,...,2\mu_n\}+\sigma>1.
\end{equation}
 Since $\min\{2\mu_1,...,2\mu_n\}>\frac 12$, such $\sigma$ can be found.
By the definition of $H_k$, (\ref{onees}) still holds.
The equation for $\hat \Phi_i^k$ becomes
$$
\Delta \hat \Phi_i^k+|y|^{2\gamma_i}e^{\xi_i^k}(\sum_j a_{ij}\hat \Phi_j^k)=\frac{O((1+|y|)^{-3-2\gamma_{n+1-i}})}{H_k(1+|y_k|)^{\sigma}},
$$
Let $y_k$ be where $H_k$ is attained. Then by the same argument as in {\bf Case one}, $|y_k|\to \infty$. In order to get a contradiction to this case, we observe that (\ref{13may14e1}) becomes
\begin{eqnarray}\label{13dec9e1}
&&\hat \Phi_i^k(y_k)-\hat \Phi_i^k(0)\\
&=&\int_{B_{\epsilon_k^{-1}}}(G_k(y_k,\eta)-G_k(0,\eta))\bigg ( |\eta |^{2\gamma_i}e^{\xi_i^k(\eta)}(\sum_j a_{ij}\tilde \Phi_j^k(\eta))
\nonumber \\
&&\qquad \qquad +\frac{O((1+|\eta |)^{-3-2\gamma_{n+1-i}})}{H_k(1+|y_k|)^{\sigma}}\bigg )d\eta
\nonumber
\end{eqnarray}
Using the same estimate on $G_k$ and (\ref{13dec9e2}) we see that the right hand side of (\ref{13dec9e1}) is $o(1)$, thus we get a contradiction as in {\bf Case one} and have proved
$$|\tilde \Phi_i^k(y)|\le C\epsilon_k(1+|y|)^{\sigma}\mbox{ for {\bf Case two} }. $$
Note that the main reason that the power of $\epsilon_k$ can be $1$ is because (\ref{13dec9e2}) holds.
Theorem \ref{thm-main} follows from the estimates of $\tilde \Phi_i^k$ and (\ref{psik}).  $\Box$

\section{The determination of $p_1,..,p_{n^2+2n}$}

In this section we explain how $p_1$, ... $p_{n^2+2n}$ are chosen to make the matrices ${\bf M}$ and ${\bf M_1}$ both invertible.

First we list some facts that can be verified easily by direct computation:
Using (\ref{nlam}) (recalling that $\tilde v^1=-\log (\lambda_0+\sum_{i=1}^n \lambda_i |P_i(z)|^2)$) we have
\begin{eqnarray}\label{cal1}
&&\frac{\partial \tilde v^1}{\partial \lambda_0}=\frac{\frac{\lambda_n}{\lambda_0}|P_n(z)|^2-1}{\lambda_0+\sum_i \lambda_i |P_i(z)|^2}, \\
&&\frac{\partial \tilde v^1}{\partial \lambda_i}=\frac{\frac{\lambda_n}{\lambda_i}|P_n(z)|^2-|P_i(z)|^2}{\lambda_0+\sum_i \lambda_i |P_i(z)|^2},\quad i=1,...,n-1, \nonumber\\
&&\frac{\partial \tilde v^1}{\partial c_{ij}^{\mathbb R}}=-\frac{2\lambda_iRe(z^{\mu_1+...+\mu_j}\bar P_i)}{\lambda_0+\sum_i \lambda_i |P_i(z)|^2}
\quad j<i, \quad i=1,..,n\nonumber \\
&&\frac{\partial \tilde v^1}{\partial c_{ij}^{\mathbb I}}=\frac{2\lambda_i Im(z^{\mu_1+...+\mu_j}\bar P_i)}{\lambda_0+\sum_i \lambda_i |P_i(z)|^2}
\quad j<i, \quad i=1,..,n\nonumber
\end{eqnarray}
It is easy to verify that for $|z|$ large
\begin{eqnarray*}
&&z^{\mu_1+..+\mu_j}\bar P_i\\
&=&|z|^{2\mu_1+..+2\mu_j+\mu_{j+1}+..+\mu_i}\bigg ( e^{-\sqrt{-1}(\mu_{j+1}+..+\mu_i)\theta}+O(|z|^{-\delta})\bigg )
\end{eqnarray*}
for some $\delta>0$ that depends only on $\mu_1,...,\mu_n$.
Thus for $|z|$ large
\begin{eqnarray}\label{cal2}
&&\qquad \frac{\partial \tilde v^1}{\partial c_{ij}^{\mathbb R}}(z) (\lambda_0+\sum_k\lambda_k|P_k(z)|^2)\\
&=&-2\lambda_i
|z|^{2\mu_1+..+2\mu_j+\mu_{j+1}+..+\mu_i}\bigg (\cos((\mu_{j+1}+..+\mu_i)\theta)+O(|z|^{-\delta}) \bigg ) \nonumber
\end{eqnarray}
\begin{eqnarray}\label{cal5}
&&\qquad \frac{\partial \tilde v^1}{\partial c_{ij}^{\mathbb I}}(z) (\lambda_0+\sum_k\lambda_k|P_k(z)|^2)\\
&=&-2\lambda_i
|z|^{2\mu_1+..+2\mu_j+\mu_{j+1}+..+\mu_i}\bigg (\sin((\mu_{j+1}+..+\mu_i)\theta)+O(|z|^{-\delta})\bigg ).\nonumber
\end{eqnarray}

By the definition of $P_i(z)$ in (\ref{piz}),
\begin{equation}\label{cal3}
|P_i(z)|^2=|z|^{2\mu_1+..+2\mu_i}(1+O(|z|^{-\delta})).
\end{equation}

We also note that
$$
\frac{\partial \tilde v^1}{\partial \lambda_i}= \frac{\lambda_0}{\lambda_i} \frac{\partial \tilde v^1}{\partial \lambda_0}
+
\frac{\frac{\lambda_0}{\lambda_i} -|P_i(z)|^2}{\lambda_0+\sum_i \lambda_i |P_i(z)|^2},\quad i=1,...,n-1.
$$

The idea of choosing $n^2+2n$ points is to make ${\bf M}$ ( ${\bf M}$ is defined in (\ref{thetamatrix})) similar to a Vandermonde type matrix. We shall use different parameters in the definition of $p_l$, which are either large or small, in order to make the leading terms dominate other terms.

Now we look at ${\bf M}$, clearly the factor $\lambda_0+\sum_k \lambda_k |P_k(p_l)|^2$ can be taken out from the $l-th$ column, thus for $|p_l|>>1$,  ${\bf M}$ is invertible if and only if
$$ {\bf M_2}:=(\Theta_1(p_1),...,\Theta_1(p_{n^2+2n})) $$
is invertible, where, according to (\ref{cal2}), (\ref{cal5}) and (\ref{cal3})
\begin{eqnarray*}
&&\Theta_1(p_l)\\
&=&\bigg (|p_l|^{2a_n}(1+O(\frac{1}{|p_l|^{\delta}}),|p_l|^{2a_{n-1}+a_{n,n-1}}\cos(a_{n,n-1}\theta_l)(1+O(\frac 1{|p_l|^{\delta}})),\\
&& |p_l|^{2a_{n-1}+a_{n,n-1}}\sin(a_{n,n-1}\theta_l)(1+O(\frac 1{|p_l|^{\delta}})),
...... \quad \bigg )'
\end{eqnarray*}
where
$$a_0=0,\quad a_i=\mu_1+..+\mu_i \,\, (i=1,...,n), \quad a_{ij}=\mu_{j+1}+..+\mu_i \,\, (i=1,..,n, j<i), $$
$\theta_l=arg(p_l)$,
$\delta>0$ only depends on $\mu_1$,..,$\mu_n$. Note that $a_{ij}=a_i-a_j$ and $ 2 a_j+a_{ij}= a_i +a_j$. The powers of $|p_l|$ are arranged in a non-decreasing order (so the largest power is $2a_n$, the second largest power is $2a_{n-1}+a_{n,n-1}$, etc). The powers of $|p_l|$ are either $2 a_i$ or $ a_i + a_j$.  Here we note that some powers appear only once (for example $2a_n$). Some powers appear only twice (for example $2a_{n-1}+a_{n,n-1}$), and it is possible that some powers appear more than twice.

Let
$$p_l=s^{1+\epsilon l} N e^{\sqrt{-1} \theta_l}, \quad l=1,..,n^2+2n $$
where $N>>s>>1>>\epsilon>0$ are constants only depending on $\mu_1,...,\mu_n$, $n$.  The angles $\theta_l$ also only depend on these parameters. We shall determine these constants and angles in the sequel.

On each row a power of $N$ can be taken out, therefore ${\bf M_2}$ is invertible iff
$$(\Theta_2(p_1),..,\Theta_2(p_{n^2+2n}))$$
is invertible, where
\begin{eqnarray*}
&&\Theta_2(p_l)=\bigg (s^{2a_n(1+\epsilon l)}(1+O(\frac{1}{|p_l|^{\delta}})),\\
&&s^{(2a_{n-1}+a_{n,n-1})(1+\epsilon l)}\cos(a_{n,n-1}\theta_l)(1+O(\frac 1{|p_l|^{\delta}})),\\
&&s^{(2a_{n-1}+a_{n,n-1})(1+\epsilon l)}\sin(a_{n,n-1}\theta_l)(1+O(\frac 1{|p_l|^{\delta}})),...,\bigg )'
\end{eqnarray*}
Hence for fixed $s$, if $N$ is sufficiently large, $O(1/|p_l|^{\delta})$ is very small, ${\bf M_2}$ is invertible iff the following matrix is invertible:
$${\bf M_3}=(\Theta_3(p_1),...,\Theta_3(p_{n^2+2n}))$$
where
\begin{eqnarray*}
&&\Theta_3(p_l)\\
&=&(s^{2a_n(1+\epsilon l)},
s^{(2a_{n-1}+a_{n,n-1})(1+\epsilon l)}\cos(a_{n,n-1}\theta_l),s^{(2a_{n-1}+a_{n,n-1})(1+\epsilon l)}\sin(a_{n,n-1}\theta_l),..)'.
\end{eqnarray*}
We start with the largest entry in ${\bf M_3}$: $s^{2a_n(1+\epsilon(n^2+2n))}$, which is in row one and column $n^2+2n$. We divide row 1 by $s^{2a_n(1+\epsilon (n^2+2n))}$ ( we call this {\bf operation one}), then the entries in row one become
$$s^{2a_n\epsilon (l-n^2-2n)}, \mbox{ for } l=1,...,n^2+2n. $$
Next we subtract a multiple of row one from other rows to eliminate the last entry in each row (we call this {\bf operation two}). For any entry in
the cofactor matrix of $1$,  if before {\bf operation two} it is of the form $s^a A$, it becomes $s^a (A+O(s^{-\delta}))$ after {\bf operation two}. Indeed, for example, let
$s^{2a_{i_0}(1+\epsilon l)}$ be an entry before {\bf operation two}. The last entry of the same row is
$s^{2a_{i_0}(1+\epsilon (n^2+2n))}$. In {\bf operation two} we subtract the $s^{2a_{i_0}(1+\epsilon (n^2+2n))}$ multiple of the first row. The entry in row 1 and the same column of $s^{2a_{i_0}(1+\epsilon l)}$ is $s^{2a_n\epsilon(l-n^2-2n)}$. Thus after {\bf operation two}
$s^{2a_{i_0}(1+\epsilon l)}$ becomes
\begin{eqnarray*}
&& s^{2a_{i_0}(1+\epsilon l)}-s^{2a_{i_0}(1+\epsilon (n^2+2n))}s^{2a_n \epsilon (l-n^2-2n)}\\
&=& s^{2a_{i_0}(1+\epsilon l)} (1-s^{(2a_{i_0}-2a_n)\epsilon (n^2+2n-l)})\\
&=& s^{2a_{i_0}(1+\epsilon l)} (1+O(s^{-\delta}))
\end{eqnarray*}
where we have used $a_{i_0}<a_n$.

Similarly if an entry before {\bf operation two} is
$$s^{(2a_j+a_{ij})(1+\epsilon l)}\cos (a_{ij}\theta_l),$$
after {\bf operation two} it becomes
$$s^{(2a_j+a_{ij})(1+\epsilon l)}(\cos (a_{ij}\theta_l)+O(s^{-\delta})), $$
for some $\delta>0$.  Eventually $s$ will be chosen large to eliminate the influence of all the perturbations.

Our strategy is to use high powers of $s$ to simplify the matrix. After the aforementioned row operations it is clear that we only need to consider the cofactor matrix of $1$, which we use $A_1$ to denote.  The highest power of $s$ in $A_1$ is shared by two entries:
$$s^{(2a_{n-1}+a_{n,n-1})(1+\epsilon(n^2+2n-1))}(\cos(a_{n,n-1}\theta_{n^2+2n-1})+O(s^{-\delta})) $$
and
$$s^{(2a_{n-1}+a_{n,n-1})(1+\epsilon(n^2+2n-1))}(\sin(a_{n,n-1}\theta_{n^2+2n-1})+O(s^{-\delta})). $$
We recall that the previous one is in row one of $A_1$.
We choose $\theta_{n^2+2n-1}=0$. In $A_1$ we divide the first row by $s^{(2a_{n-1}+a_{n,n-1})(1+\epsilon (n^2+2n-1))}$, then the largest entry in
row $1$ of $A_1$ becomes $1+O(s^{-\delta})$. We then subtract from other rows a multiple of the first row to eliminate the last entry of each row.  By the same reason as before, after these row operations the invertibility of $A_1$ is equivalent to the invertibility of the cofactor matrix $A_2$ of
$1+O(s^{-\delta})$, a  $(n^2+2n-2)\times (n^2+2n-2)$ matrix which is barely changed after these transformations. In fact, each entry in $A_2$ is only multiplied a factor $1+O(s^{-\delta})$ after these transformations.

As we continue this process we face three situations. If the highest power of $s$ without the $\epsilon$ part is not repeated, we just apply the same type of row operations as in {\bf operation one} and {\bf operation two}. If the highest power of $s$ without the $\epsilon$ part is shared by only
two entries (one is a cosine term, one is a sine term), we just take the corresponding angle to be $0$, so the cosine term will dominate all other terms and this case is similar to the previous case. Finally we may run into the following situation: A power of $s$ without the $\epsilon$ part is shared by more than two indices:
\begin{eqnarray*}
 &&\exists i_0,j_0,i_1,j_1, \mbox{ such that } 2a_{j_0}+a_{i_0,j_0}=2a_{j_1}+a_{i_1,j_1},\quad j_0\neq j_1. \\
 && \exists i_0,j_0,i_1, \mbox{ such that } 2a_{j_0}+a_{i_0,j_0}=2a_{j_1}.
 \end{eqnarray*}

 In this case we first prove the following simple but important lemma.

 \begin{lem}\label{cal-lem1} There exist $\epsilon_0>0$ that depends only on $\mu_1,..,\mu_n$ and $n$ such that
 for $\epsilon\in (0,\epsilon_0)$,
 \begin{equation}\label{fact1}
 \frac{|p_a|^{l_1}}{|p_b|^{l_2}}\to \infty \mbox{ as } s\to \infty, \forall a, b \in \{1,...,n^2+2n\}.
 \end{equation}
 where $l_1$, $l_2$ are two numbers in the set $\{ 2a_1,..,2a_n, ..., 2a_j+ a_{ij},...,\}$ that satisfy $l_1>l_2$.
 \end{lem}

 \noindent{\bf Proof of Lemma \ref{cal-lem1}:}
 Suppose
 $|p_a|^{l_1}=s^{(1+\epsilon a)l_1}$, $|p_b|^{l_2}=s^{(1+\epsilon b)l_2}$, it is easy to see that for all $a,b\in \{1,..,n^2+2n\}$,
 $(1+\epsilon a)l_1>(1+\epsilon b)l_2$ if $l_1>l_2$ and $\epsilon$ is sufficiently small. The smallness of $\epsilon$ is clearly determined by
 the set
 $$\{2a_1,..,2a_n, ..., 2a_j+ a_{ij},...,\}. $$
 Lemma \ref{cal-lem1} is established. $\Box$

 \medskip

Next we prove two more Calculus lemmas.

\begin{lem}\label{cal-lem}
Let $N_1<N_2<...<N_k$ be positive numbers. Then there exist $\theta_1,\theta_2,...,\theta_{2k+1}$ such that the following matrix
$$M_{Nk}=\left (\begin{array}{cccccc}
1 & ... & ... & ... & ... & 1 \\
\sin(N_1\theta_1) & ... & ... & ... & ... & \sin(N_1\theta_{2k+1}) \\
\cos(N_1\theta_1) & ... & ... & ... & ... & \cos(N_1\theta_{2k+1}) \\

... & ... & ... & ... & ... & ... \\
\sin(N_k\theta_1) & ... & ... & ... & ... & \sin(N_k\theta_{2k+1})\\
\cos(N_k\theta_1) & ... & ... & ... & ... & \cos(N_k\theta_{2k+1})
\end{array}
\right )
$$
satisfies
$$0<c_1(N_1,...,N_k)<| det(M_{Nk}) |<c_2(N_1,...,N_k). $$
for positive constants $c_1$ and $c_2$ that only depend on $N_1,...,N_k$.
\end{lem}

\noindent{\bf Proof of Lemma \ref{cal-lem}:} We use the Taylor expansion of $\sin(N\theta)$ and $\cos(N\theta)$:
$$\sin (N_i\theta_j)=\sum_{l=1}^k (-1)^{l+1}\frac{(N_i \theta_j)^{2l-1}}{(2l-1)!}+O((N_i \theta_j )^{2k+1}). $$
$$\cos (N_i\theta_j)=\sum_{l=0}^k (-1)^{l}\frac{(N_i \theta_j)^{2l}}{(2l)!}+O((N_i \theta_j )^{2k+2}). $$
We apply the following elementary operations on $M_{Nk}$: First we subtract a multiple of the first row from other odd number rows to eliminate
the first order terms of $\theta_i$ ($i=1,..,2k+1$). After the cancelation it is easy to see that, the entry of row $2j-1$ ($j>1$) and column $r$ ($r>1$) is of the form
$$\sum_{l=2}^k (-1)^{l+1}(a_{l,j}\theta_r)^{2l-1}+O(\theta_r)^{2k+1}$$
for some positive constant $a_{l,j}$, which satisfies
$a_{l,j}<a_{l,j+1}$. In the second step we use row three to eliminate all the $O(\theta^3)$ terms of other odd number rows starting from row 5.
After the second step, the entry of row $2j-1$ ($j>2$) and column $r$ ($r>2$) is of the form
$$
\sum_{l=3}^k (-1)^{l+1}(\tilde a_{l,j}\theta_r)^{2l-1}+O(\theta_r)^{2k+1},
$$
with $\tilde a_{l,j}>0$ satisfying $\tilde a_{l,j}<\tilde a_{l,j+1}$.

After $k-1$ such operations we see that the entry of
row $2j-1$ and column $r$ is a multiple of $\theta_r^{2j-1}$ plus lower order terms. Clearly we can use the terms on row $2k-1$ to eliminated all the $O(\theta^{2k-1})$ terms in other odd number rows. Then we can use row $2k-3$ to remove the $O(\theta^{2k-3})$ terms in other odd number rows. After such operations the entry of row $2j-1$ and column $r$ is
$ C\theta_r^{2j-1}+O(\theta_r^{2k+1})$. Similar operations can be applied to even number rows. Thus after a finite number of elementary row operations ( including multiplying a constant on each row) the matrix $M_{Nk}$ is transformed to
$$\tilde M_{Nk}=\left (\begin{array}{cccccc}
1 & 1 & ... & ... & ... & 1 \\
\theta_1 & \theta_2 & ... & ... & ... & \theta_{2k+1} \\
\theta_1^2 & \theta_2^2 & ... & ... & ... & \theta_{2k+1}^2 \\
... & ... & ... & ... & ... & ... \\
\theta_1^{2k-1} & \theta_2^{2k-1} & ... & ... & ... & \theta_{2k+1}^{2k-1}\\
\theta_1^{2k} & \theta_2^{2k} & ... & ... & ... & \theta_{2k+1}^{2k}
\end{array}
\right )  + \mbox{ a minor matrix }.
$$
The $(i,j)$ entry of the second matrix is $O(\theta_i^{2k+1})$. Now we choose $\theta_i=i\epsilon$ for some $\epsilon>0$ that depends only on $N_1,...,N_k$. For $\epsilon$ sufficiently small, $\tilde M_{Nk}$ is invertible if and only if the first matrix is invertible. Finally we observe that the first matrix of $\tilde M_{Nk}$ is a Vandermonde matrix.
Lemma \ref{cal-lem} is established. $\Box$

\medskip

The proof of the following lemma is very similar and is omitted.

\begin{lem}\label{cal-lem2}
Let $N_1<N_2<...<N_k$ be positive numbers. Then there exist $\theta_1,\theta_2,...,\theta_{2k}$ such that the following matrix
$$ M_{2Nk}=\left (\begin{array}{cccccc}
\sin(N_1\theta_1) & ... & ... & ... & ... & \sin(N_1\theta_{2k}) \\
\cos(N_1\theta_1) & ... & ... & ... & ... & \cos(N_1\theta_{2k}) \\

... & ... & ... & ... & ... & ... \\
\sin(N_k\theta_1) & ... & ... & ... & ... & \sin(N_k\theta_{2k})\\
\cos(N_k\theta_1) & ... & ... & ... & ... & \cos(N_k\theta_{2k})
\end{array}
\right )
$$
satisfies
$$0<c_1(N_1,...,N_k)<| det( M_{2Nk}) |<c_2(N_1,...,N_k). $$
for positive constants $c_1$ and $c_2$ that only depend on $N_1,...,N_k$.
\end{lem}

\bigskip

Now we go back to the case that after finite steps of reduction, the highest power of $s$ without the $\epsilon$ part is $M$ and is shared by more than $2$ indices. Our goal is to make the following matrix invertible:

$${\bf A_2}=\left (\begin{array}{cc}
B & C \\
D & F
\end{array}
\right )  \cdot (1+O(s^{-d}))
$$
where the last term $(1+O(s^{-d}))$ means each entry in $\displaystyle{\left (\begin{array}{cc}
B & C \\
D & F
\end{array}
\right )} $ is multiplied by a quantity of the magnitude $1+O(s^{-d})$, even though
these quantities are different from one another.
$C$ is either of the form
$$\left(\begin{array}{ccc}
s^{M(1+\epsilon (l+1))}\sin (N_1\theta_{l+1}) & ... & s^{M(1+\epsilon (l+2T))}\sin (N_1\theta_{l+2T}) \\
s^{M(1+\epsilon (l+1))}\cos (N_1\theta_{l+1}) & ... & s^{M(1+\epsilon (l+2T))}\cos (N_1\theta_{l+2T}) \\
... & ... & ... \\
s^{M(1+\epsilon (l+1))}\sin (N_T\theta_{l+1}) & ... & s^{M(1+\epsilon (l+2T))}\sin (N_T\theta_{l+2T})\\
s^{M(1+\epsilon (l+1))}\cos (N_T\theta_{l+1}) & ... & s^{M(1+\epsilon (l+2T))}\cos (N_T\theta_{l+2T})
\end{array}
\right )
$$
or
$$\left(\begin{array}{ccc}
1 & ... & 1\\
s^{M(1+\epsilon (l+1))}\sin (N_1\theta_{l+1}) & ... & s^{M(1+\epsilon (l+2T+1))}\sin (N_1\theta_{l+2T+1}) \\
s^{M(1+\epsilon (l+1))}\cos (N_1\theta_{l+1}) & ... & s^{M(1+\epsilon (l+2T+1))}\cos (N_1\theta_{l+2T+1}) \\
... & ... & ... \\
s^{M(1+\epsilon (l+1))}\sin (N_T\theta_{l+1}) & ... & s^{M(1+\epsilon (l+2T+1))}\sin (N_T\theta_{l+2T+1})\\
s^{M(1+\epsilon (l+1))}\cos (N_T\theta_{l+1}) & ... & s^{M(1+\epsilon (l+2T+1))}\cos (N_T\theta_{l+2T+1})
\end{array}
\right )
$$
We take the first case as an example. $B$ is of the form
$$B=\left(\begin{array}{ccc}
s^{M(1+\epsilon)}\sin (N_1\theta_1) & ... & s^{M(1+\epsilon l)}\sin (N_1\theta_l) \\
s^{M(1+\epsilon)}\cos (N_1\theta_1) & ... & s^{M(1+\epsilon l)}\cos (N_1\theta_l) \\
... & ... & ... \\
s^{M(1+\epsilon )}\sin (N_T\theta_1) & ... & s^{M(1+\epsilon l)}\sin (N_T\theta_l)\\
s^{M(1+\epsilon)}\cos (N_T\theta_1) & ... & s^{M(1+\epsilon l)}\cos (N_T\theta_l)
\end{array}
\right )
$$
The importance of Lemma \ref{cal-lem1} is that it makes $F$ minor.
For matrices $D$ and $F$, we just write one row vector of $(D, F)$ as a representative:
$$ \bigg (s^{H(1+\epsilon)}, ..., s^{H(1+\epsilon l)}, s^{H(1+\epsilon(l+1))},..., s^{H(1+\epsilon (l+2T))} \bigg ) $$
where
$$\bigg (s^{H(1+\epsilon)}, ..., s^{H(1+\epsilon l)} \bigg )$$
is a row vector of $D$,
$$\bigg (s^{H(1+\epsilon(l+1))},..., s^{H(1+\epsilon (l+2T))} \bigg ) $$
is a row vector of $F$.
Here we note that $H<M$, other rows of ${\bf A_2}$ may have sine or cosine terms.

Now we take $s^{M(1+\epsilon (l+1))}$ out of the $2k$ rows of $(B, C)$, after this operation $B$ and $C$ become $\tilde B$ and $\tilde C$:
$$
\tilde B=\left(\begin{array}{cccc}
s^{-M\epsilon l}\sin (N_1\theta_l) & s^{-M\epsilon (l-1)}\sin (N_1\theta_2) & ... & s^{-M\epsilon }\sin (N_1\theta_l) \\
s^{-M\epsilon l}\cos (N_1\theta_1) & s^{-M\epsilon (l-1)}\cos (N_1\theta_2) & ... & s^{-M\epsilon }\cos (N_1\theta_l) \\
... & ... & ... & ...  \\
s^{-M\epsilon l}\sin (N_T\theta_l) & s^{-M\epsilon (l-1)}\sin (N_T\theta_2) & ... & s^{-M\epsilon }\sin (N_T\theta_l)\\
s^{-M\epsilon l}\cos (N_T\theta_1) & s^{-M\epsilon (l-1)}\cos (N_T\theta_2) & ... & s^{-M\epsilon }\cos (N_T\theta_l)
\end{array}
\right )
$$
$$
\tilde C=\left(\begin{array}{cccc}
\sin (N_1\theta_{l+1}) & s^{M\epsilon}\sin (N_1\theta_{l+2}) & ... & s^{M(2T-1)\epsilon}\sin (N_1\theta_{l+2T}) \\
\cos (N_1\theta_{l+1}) & s^{M\epsilon}\cos (N_1\theta_{l+2}) & ... & s^{M(2T-1)\epsilon}\cos (N_1\theta_{l+2T}) \\
... & ... & ... & ... \\
\sin (N_T\theta_{l+1}) & s^{M\epsilon}\sin (N_T\theta_{l+2}) & ... & s^{M(2T-1)\epsilon}\sin (N_T\theta_{l+2T})\\
\cos (N_T\theta_{l+1}) & s^{M\epsilon}\cos (N_T\theta_{l+2}) & ... & s^{M(2T-1)\epsilon}\cos (N_T\theta_{l+2T})
\end{array}
\right )
$$
After these row operations the major part of ${\bf A_2}$ becomes
$${\bf A_3}=(A_{31}, A_{32})=\left ( \begin{array}{cc}
\tilde B & \tilde C \\
D & F
\end{array}
\right )
$$
Starting from the second column of $A_{32}$ we take away the power of $s$. For example we divide the second column of $A_{32}$ by $s^{M\epsilon}$,
the third column by $s^{2M\epsilon }$ and the $2T-th$ column by $s^{M(2T-1)\epsilon}$.
Now we see the influence of the representative row vector in $F$. Before this set of column operations it is
$$\bigg (s^{H(1+\epsilon (l+1))}, ..., s^{H(1+\epsilon (l+2T))} \bigg ) $$
After these column operations it becomes (using $H<M$)
$$ s^{H(1+\epsilon (l+1))} \bigg (1, O(s^{-d}),..., O(s^{-d}) \bigg ). $$
Note that this computation is very similar to those in the proof of Lemma \ref{cal-lem1}.
We use $\tilde F$ to represent the new matrix after the column operations on $F$.

After these column operations, $\tilde C$ becomes
$$\tilde C_1=\left(\begin{array}{cccc}
\sin (N_1\theta_{l+1}) & \sin (N_1\theta_{l+2}) & ... & \sin (N_1\theta_{l+2T}) \\
\cos (N_1\theta_{l+1}) & \cos (N_1\theta_{l+2}) & ... & \cos (N_1\theta_{l+2T}) \\
... & ... & ... & ... \\
\sin (N_T\theta_{l+1}) & \sin (N_T\theta_{l+2}) & ... & \sin (N_T\theta_{l+2T})\\
\cos (N_T\theta_{l+1}) & \cos (N_T\theta_{l+2}) & ... & \cos (N_T\theta_{l+2T})
\end{array}
\right )
$$
By Lemma \ref{cal-lem2}, $\tilde C_1$ is invertible, which means its row vectors are linearly independent. Thus there is a combination of its row vectors to cancel the representative vector in $\tilde F$ (just the major part):
$$s^{H(1+\epsilon (l+1))} \bigg (1, 0,..., 0 \bigg ). $$
When this same row operation is applied to $A_{31}$, the representative vector in $D$:
$$(s^{H(1+\epsilon)},..., s^{H(1+\epsilon l)}) $$
becomes this after the row transformation:
$$(s^{H(1+\epsilon)}(1+O(s^{-d})),..., s^{H(1+\epsilon l)}(1+O(s^{-d}))) $$
where we used $H<M$ again. After these elementary operations, $B$ and $F$ are turned into minor matrices. Thus the invertibility of ${\bf A_2}$ is reduced to the invertibility of the transformation of $D$, which is of the same nature of $D$. This method of reduction can be continued and the construction of $p_1$, ... $p_{n^2+2n}$ is complete for matrix ${\bf M}$.

Since ${\bf M_1}$ is very similar to ${\bf M}$ and we only require $N$, $s$ to be large and $\epsilon$ to be small in ${\bf M_1}$. Moreover the angles in ${\bf M_1}$ are the same as in ${\bf M}$. Thus $p_1,..., p_{n^2+2n}$ that make ${\bf M}$ invertible also make ${\bf M_1}$ invertible.
The construction of $p_1,.., p_{n^2+2n}$ is complete.

\bigskip
\footnotesize
\noindent\textit{Acknowledgments.}
The research of Wei is supported by a GRF from RGC of Hong Kong and NSERC of Canada. Part of the paper was finished when the third author was visiting Chinese University of Hong Kong from April to May in 2012 and
Taida Institute of Mathematical Sciences (TIMS) in June 2012 and July 2013. He would like to thank both institutes for their warm hospitality and generous financial support.

\end{document}